\documentclass[12pt]{article}

\providecommand{\keywords}[1]{\textbf{\textit{Keywords:}} #1}

\usepackage{amsthm}
\usepackage[english]{babel}
\usepackage{amsthm}

\newtheorem{theorem}{Theorem}[section]
\newtheorem{lemma}[theorem]{Lemma}

\usepackage{amsmath}
\usepackage{subfig}
\usepackage{graphicx}
\usepackage{float}
\usepackage{times}
\usepackage{bm}
\usepackage{natbib}
\usepackage[justification=centering]{caption}
\usepackage[plain,noend]{algorithm2e}

\usepackage{times}
\usepackage{bm}
\usepackage{natbib}
\usepackage[justification=centering]{caption}
\usepackage[plain,noend]{algorithm2e}
\newcommand{\udots}{\mathinner{\mskip1mu\raise1pt\vbox{\kern7pt\hbox{.}}  
		\mskip2mu\raise4pt\hbox{.}\mskip2mu\raise7pt\hbox{.}\mskip1mu}}    
\makeatletter
\renewcommand{\algocf@captiontext}[2]{#1\algocf@typo. \AlCapFnt{}#2} % text of caption
\def\@algocf@capt@plain{top}
\renewcommand{\algocf@makecaption}[2]{%
  \addtolength{\hsize}{\algomargin}%
  \sbox\@tempboxa{\algocf@captiontext{#1}{#2}}%
  \ifdim\wd\@tempboxa >\hsize%     % if caption is longer than a line
    \hskip .5\algomargin%
    \parbox[t]{\hsize}{\algocf@captiontext{#1}{#2}}% then caption is not centered
  \else%
    \global\@minipagefalse%
    \hbox to\hsize{\box\@tempboxa}% else caption is centered
  \fi%
  \addtolength{\hsize}{-\algomargin}%
}
\makeatother

\begin{document}

%% Here are the title, author names and addresses
\title{Te Test: A New Non-asymptotic T-test for Behrens-Fisher Problem}

\author{Chang Wang, Jinzhu Jia}
\date{}
\maketitle

\begin{abstract}

The Behrens-Fisher Problem is a classical statistical problem. It is to test the equality of the means of two normal populations using two independent samples, when the equality of the population variances is unknown. \cite{linnik2008} has shown that this problem has no exact fixed-level tests based on the complete sufficient statistics. However, exact conventional solutions based on other statistics and approximate solutions based the complete sufficient statistics do exist.\citep{Weerahandi1995Exact} 

Existing methods are mainly asymptotic tests, and usually don't perform well when the variances or sample sizes differ a lot. In this paper, we propose a new method to find an exact t-test ($T_e$) to solve this classical Behrens-Fisher Problem. Confidence intervals for the difference between two means are provided. We also use detailed analysis to show that $T_e$ test reaches the maximum of degree of freedom and to give a weak version of proof that $T_e$ test has the shortest confidence interval length expectation. Some simulations are performed to show the advantages of our new proposed method compared to available conventional methods like Welch's test, paired t-test and so on. We will also compare it to unconventional method, like two-stage test.
\end{abstract}


\section{Introduction}

Test of the equality of the means of  two normal populations is a classical statistical problem. It is well known and well accepted  when the variances of the two populations are the same but unknown, a t-test could be used. But when the variances are not the same, which is very likely the case in practice, this problem is not solved perfectly. Although a few text books suggest testing the equality of variances first, it is still problematic even if we could not reject the null hypothesis that the variances are equal, because it doesn't necessarily mean the equality of variance. The test of the equality of two means when the variances are not the same is called the {\it{Behrens-Fisher Problem}} \citep{welch1938significance}. A lot of researchers have studied this problem or improve the solutions \citep{Best1987Welch,Fenstad1983A,twostage2,Paul1992COMMENT,comparison1}. 

When the sample sizes in each of the two population are big, it is not hard to understand that many asymptotic methods perform well. For example, the standard
normal statistics ($T_N$) are recommended when sample sizes are big; when the sample sizes are small and close, the two sample t-test ($T_1$) with Welch's degree of freedom performs well \citep{comparison1,ReviewT1,welch}. %A bootstrap method for this problem is also developed by (Efron and Tibshirani, 1993)[13]. 
%Besides, although non-parametric and Bayesian procedures are also being developed recently, many evidences are in favor of Welch's $T_1$. It is alos highly recommanded by some text books (Hogg and Tanis, 2010).
A few other methods are introduced and compared in a review paper \citep{ReviewT1}. Different from the previous test methods, most of which are approximate tests, our method is an exact $t$ test with well defined degree of freedom. We call our method $T_e$ (e for exact) test.

This article is organized as follows. In section 2, we introduce some methods available now, and we will compare them with our proposed method later.
 In section 3, we introduce the formula of our new method and prove its usefulness. 
 In section 4 and 5, we use some important theorems to show the theoretical superiority of $T_e$ test. 
 In section 6, we use simulations to show the advantages and disadvantages of our method. 
 In section 7, we develop the $T_e$ test to the condition of multi-variate form by introducing Hotelling's T square.
 What's more, we also provide an unconventional exact test of Behrens Fisher Problem for comparison in section 8.
 Finally, we conclude in section 9.

\section{Methods available}
%\subsection{problem introduction}
Suppose now we have two independent normal samples $X_i, i \in \{1,2,\ldots,m\}$ and $Y_j, j \in \{ 1,2,\ldots,n\}$, where  $X_i \sim N(\mu_1,\sigma^2_1), Y_j \sim N(\mu_2,\sigma^2_2)$. Our objective is to find a test method to test the null hypothesis that they have the same expectation ($\mu_1-\mu_2 = 0$).   In the following discussions, we  assume that $m\geq n$.  A few useful statistics are defined below. \begin{equation*}
\bar X = m^{-1} \sum_{i=1}^m X_j,\quad s_1^2 = \sum_{i=1}^m (X_i-\bar X)^2;
\end{equation*}
\begin{equation*}
\bar Y = n^{-1} \sum_{j=1}^n Y_j,\quad s_2^2 = \sum_{j=1}^n (Y_j-\bar Y)^2;
\end{equation*}

Now we briefly introduce a few commonly used test statistics for this classic problem\citep{ReviewT1}. They are $T_N$ test, Welch's $T_1$ test and paired $t$ test.
%\subsection{$T_{N}$ test}

\subsection{$T_N$ test}
When both $n_1$ and $n_2$ are sufficiently large,   $T_N$ defined below
\begin{equation*}
T_N=\frac{\bar{X}-\bar{Y}-(\mu_1-\mu_2)}{\sqrt{\frac{s_1^2}{m}+\frac{s_2^2}{n}}}
\end{equation*}
 is asymptotically normal with mean 0 and variance 1.  Test using an approximate normal statistics $T_N$ is called $T_N$ test. Apparently, it is useful only when both $m$ and $n$ are big.

\subsection{$T_1$ test}
 
When the sample sizes are small, it is shown that $T_N$ is approximately distributed as students' $t$  with degree of freedom $f$ defined below:
\begin{equation*}
f=\frac{      (\frac{s_1^2}{m}+\frac{s_2^2}{n})^2       }{   \frac{s_1^4}{m^2(m-1)}+\frac{s_2^4}{n^2(n-1)}         }.
\end{equation*}
 Test using this $t$ statistics is called $T_1$ test (also called Welch's t-test) and is regarded as one of the best methods to date.

Its idea can be summarized as follow:

Let $$\sigma_B^2 = \frac{\sigma_1^2}{n_1} +  \frac{\sigma_2^2}{n_2} \quad 
s_B^2 = \frac{s_1^2}{n_1} +  \frac{s_2^2}{n_2} $$

We hope to get a $f$ that satisfy the approximation $s_B^2/\sigma_B^2 \sim \chi_f^2/f$:
\begin{align*}
	&var(\frac{rs_B^2}{\sigma_B^2})=2f\\
	\Rightarrow &\frac{f^2}{\sigma_B^4}var(\frac{s_1^2}{n_1} +  \frac{s_2^2}{n_2} )=2f\\
	\Rightarrow & f^2 ( \frac{2}{n_1^2} \frac{\sigma_1^4}{n_1-1}  + \frac{2}{n_2^2} \frac{\sigma_2^4}{n_2-1} ) = 2f\sigma_B^4\\
	\Rightarrow & f=\frac{\sigma_B^4}{\frac{1}{n_1-1} (\frac{\sigma_1^2}{n_1})^2  + \frac{1}{n_2-1} (\frac{\sigma_2^2}{n_2})^2 } \approx 
	\frac{s_B^4}{\frac{1}{n_1-1} (\frac{s_1^2}{n_1})^2  + \frac{1}{n_2-1} (\frac{s_2^2}{n_2})^2 } 
\end{align*}

which explains $T_1$ test's degree of freedom $f$.
\subsection{Paired t-test}

Although the paired t-test requires that the data are paired,  it still values as an exact method when it comes to this problem. Even if $m>n$, we can select $n$ data points from $X$ randomly to do a paired $t$ test with $Y$, it fails to use all information of $X$ though. We also use this as a benchmark.%, and we will show later that when $m=n$, the $T_e$ test is just the paired $t$-test. It falls short of our expectation because when the sample sizes are seriously unbalanced, it can only use a small number of imformation from two populations. For example, when (m,n)=(10,100), it can only use (10+10)=20 samples from 110 samples.

Let $c_i=X_i-Y_i \sim N(\mu_1-\mu_2,\sigma_1^2+\sigma_2^2)$, then:
\begin{equation*}
\frac{ \frac{1}{n}\sum_{i=1}^{n}c_i-(\mu_1-\mu_2)}    {\sqrt{\frac{\sum_{i=1}^{n}(c_i-\bar{c})^2}{n(n-1)}}}\sim t_{n-1 }
\end{equation*}

\subsection{Other tests}
There are some other conventional methods like likelihood ratio test, bootstrap methods \citep{Efron1993An}, non-parametric methods \citep{Fligner1981Robust}. But according to \cite{ReviewT1}, to the date, the statistic $T_1$ is the best and is referred as the statistic to use in recent text books. Thus we are not focus on these methods and they are omitted here. 

Besides, there are also some other unconventional methods like two-stage test and generalized pivotal. We will introduce their properties in detail later. Generalized pivotal \citep{Weerahandi1995Exact}'s performance will also be included in simulations.

We now turn to our proposed new exact test: $T_e$ test.
\section{$T_e$ test}

\subsection{Formula of $T_e$ test}
%Let $(Q^T)_{n\times m}$ be the first $n$ rows of an $m\times m$ orthogonal matrix.
%\begin{equation*}       
%(Q^T)_{n\times m}=
%\left(                 
%\begin{array}{cccccc}   
%\frac{1}{\sqrt{m}}			&\cdots&\frac{1}{\sqrt{m}}&\frac{1}{\sqrt{m}}&\frac{1}{\sqrt{m}}	 &\frac{1}{\sqrt{m}}		\\
%\frac{1}{\sqrt{m(m-1)}}		&\cdots&\frac{1}{\sqrt{m(m-1)}}&\frac{1}{\sqrt{m(m-1)}}&\frac{1}{\sqrt{m(m-1)}}&\frac{-(m-1)}{\sqrt{m(m-1)}}\\
%\frac{1}{\sqrt{(m-1)(m-2)}}	&\cdots&\frac{1}{\sqrt{(m-1)(m-2)}}&\frac{1}{\sqrt{(m-1)(m-2)}}&\frac{-(m-2)}{\sqrt{(m-1)(m-2)}}&0\\
%\vdots	&\vdots&\vdots &\udots&\vdots&\vdots\\
%\frac{1}{\sqrt{(m-n+2)(m-n+1)}}	& \cdots &\frac{-(m-n+1)}{\sqrt{(m-n+2)(m-n+1)}} &0&\cdots&0\\
%\end{array}
%\right)                 
%\end{equation*}

Let $(P^T)_{n\times n}$ be an $n\times n$ orthogonal matrix whose elements of the first row are all $1/\sqrt{n}$, and an example is:
\begin{equation*}    
(P^T)_{n\times n}=  
\left(                 
\begin{array}{cccccc}  
\frac{1}{\sqrt{n}}			&\cdots&\frac{1}{\sqrt{n}}&\frac{1}{\sqrt{n}}&\frac{1}{\sqrt{n}}&\frac{1}{\sqrt{n}}		\\
\frac{1}{\sqrt{n(n-1)}}		&\cdots&\frac{1}{\sqrt{n(n-1)}}&\frac{1}{\sqrt{n(n-1)}}&\frac{1}{\sqrt{n(n-1)}}&\frac{-(n-1)}{\sqrt{n(n-1)}}\\
\frac{1}{\sqrt{(n-1)(n-2)}}	&\cdots&\frac{1}{\sqrt{(n-1)(n-2)}}&\frac{1}{\sqrt{(n-1)(n-2)}}&\frac{-(n-2)}{\sqrt{(n-1)(n-2)}}&0\\
\vdots	&\vdots&\vdots &\udots&\vdots&\vdots\\
\frac{1}{\sqrt{2\times 3}}	 &\frac{1}{\sqrt{2\times 3}} &\frac{-2}{\sqrt{2\times 3}} &0&\cdots&0\\
\frac{1}{\sqrt{1\times 2}}	 &\frac{-1}{\sqrt{1\times 2}} &0&0&\cdots&0\\
\end{array}
\right)                 
\end{equation*}

Similarly, let $(Q^T)_{n\times m}$ be the first $n$ rows of an $m\times m$ orthogonal matrix (whose elements of the first row are all $1/\sqrt{m}$).

Let $X = [X_1,X_2,\ldots,X_m]^T$ and $Y = [Y_1,Y_2,\ldots,Y_n]^T$ be the sample vectors defined in section 2.
Then $Z:=(Q^T)_{n\times m}X/\sqrt{m}-(P^T)_{n\times n}Y/\sqrt{n}$ is an $n$-dimensional normal random vector.
And it is easy to get its expectation and variance-covariance matrix:
\begin{align*}
E(Z)&=\mu_1Q^T1_m/\sqrt{m}-\mu_2P^T1_n/\sqrt{n}=(\mu_1-\mu_2,0,...,0)^T
\end{align*}
\begin{align*}
var(Z)&=var(Q^TX/\sqrt{m}-P^TY/\sqrt{n})\\
&=var(Q^TX/\sqrt{m})+var(P^TY/\sqrt{n})\\
&=Q^Tvar(X)Q/m+P^Tvar(Y)P/n\\
&=(\frac{\sigma_1^2}{m}+\frac{\sigma_2^2}{n})I_n
\end{align*}

Let $Z = [Z_1,Z_2,\ldots,Z_n]^T$.  From the above calculation, we can get the distribution of $Z$:

$$
Z
\sim N(
(\mu_1-\mu_2,0,...,0)^T
, (\frac{\sigma_1^2}{m}+\frac{\sigma_2^2}{n})I_n). $$
%$$\begin{pmatrix}
%Z_1\\
%Z_2\\
%\vdots\\
%Z_n
%\end{pmatrix}
%\sim N(
%\begin{pmatrix}
%\mu_1-\mu_2\\
%0\\
%\vdots\\
%0
%\end{pmatrix}, (\frac{\sigma_1^2}{m}+\frac{\sigma_2^2}{n})I_n). $$

From the above distribution we see that 

$$ Z_1-(\mu_1-\mu_2)\sim N(0,\frac{\sigma_1^2}{m}+\frac{\sigma_2^2}{n}),\   \frac{\sum_{i=2}^n Z^2_i}{n-1}\sim \frac{\chi^2_{n-1}}{n-1}\times(\frac{\sigma_1^2}{m}+\frac{\sigma_2^2}{n})$$ 

 $$Z_1-(\mu_1-\mu_2)\quad  \bot \quad \sum_{i=2}^n Z^2_i$$.

So,
$$T_e := \frac{    Z_1-(\mu_1-\mu_2)        }{ \sqrt{ \sum_{i=2}^{n} Z^2_i / (n-1) }    } \sim t_{n-1}.$$

\subsection{Comparison with paired t-test}

When the sample sizes are the same, i.e, $m=n$, we can simplify the $T_e$ test statistics. At this situation,
\begin{equation*}
Z=(P^T)_{n\times n}(X-Y)/\sqrt{n}. 
\end{equation*}

Let $c_i=X_i-Y_i \sim N(\mu_1-\mu_2,\sigma_1^2+\sigma_2^2)$, then 

$$\begin{pmatrix}
Z_1\\
Z_n\\
Z_{n-1}\\
\vdots\\
Z_2\\
\end{pmatrix}=
\begin{pmatrix}
\sum_{i=1}^n c_i/n\\
(c_1-c_2)/(\sqrt{1\times 2\times n})\\
(c_1+c_2-2c_3)/(\sqrt{2\times 3\times n}) \\
\vdots\\
(\sum_{i=1}^{n-1} c_i - (n-1)c_n)/(\sqrt{(n-1)\times n\times n})
\end{pmatrix}
\sim N(
\begin{pmatrix}
\mu_1-\mu_2\\
0\\
0\\
\vdots\\
0\\
\end{pmatrix},  
\frac{\sigma_1^2+\sigma_2^2}{n}I_n
) $$
$$Z_1=\frac{\sum_{i=1}^{n}c_i}{n}$$
\begin{align*}
	\sum_{i=2}^{n} Z^2_i&=\frac{\sum_{i=1}^{n} (1-\frac{1}{n})c_i^2 -\sum_{i\not = j}c_ic_j}{n}\\
	&=\frac{\sum_{i=1}^{n} c_i^2 - \frac{1}{n}(\sum_{i=1}^{n} c_i)^2 }{n}\\
	&= \frac{\sum_{i=1}^{n}(c_i-\bar{c})^2}{n}     
\end{align*}

$$ \frac{    Z_1-(\mu_1-\mu_2)        }{ \sqrt{ \sum_{i=2}^{n} Z^2_i / (n-1) }    }=\frac{ \frac{1}{n}\sum_{i=1}^{n}c_i-(\mu_1-\mu_2)}    {\sqrt{ \sum_{i=1}^{n}(c_i-\bar{c})^2/n(n-1) }}\sim t_{n-1}    $$
which is exactly the paired t-test!

Besides, when sample sizes are not equal, although the degree of freedom of $T_e$ test is decided by the smaller sample size and is thus same as paired t-test, the variance of $T_e$ test's transformed data is smaller than that of paired t-test.

For example, suppose we have two independent normal samples $X_i, i=1,2,\ldots,m$ and $Y_j,j = 1,2,\ldots,n$ (and we assume $m\ge n$).  $$X_i \sim N(\mu_1,\sigma^2_1)\quad (i=1,...,m)\quad  Y_j \sim N(\mu_2,\sigma^2_2)\quad (j=1,...,n)$$

Before we perform the final t-test, the paired one's transformation is 
$$c_i=X_i-Y_i \sim N(\mu_1-\mu_2,\sigma_1^2+\sigma_2^2)$$
$$\frac{\sum_{i=1}^{n}c_i}{n}-(\mu_1-\mu_2) \sim N(0,\frac{\sigma_1^2+\sigma_2^2}{n})$$

While $T_e$ test's transformation is 
$$ Z_1-(\mu_1-\mu_2)\sim N(0,\frac{\sigma_1^2}{m}+\frac{\sigma_2^2}{n}),$$ 

Therefore, after transformation, $T_e$ test's data variance is smaller than that of paired t-test, which will lead to a greater power in the final t-test.

\section{Maximum of degree of freedom}

In $T_e$ test, we find that test's degree of freedom is decided by the smaller sample size. In fact, \cite{Scheff1943On} proved an important theorem, claiming that with certain form, the maximum number of degree of freedom is no greater than $min(m-1,n-1)$, where m and n is the sample sizes of the two groups. He also gave a special and simple transformation to reach it, without showing a clear division of data information, though.

\begin{theorem}[Scheffe's]\label{scheffe}
	Let L be a linear form and Q a quadratic form in the variates $x_1$, ..., $x_m$, $y_1$, ..., $y_n$ with coefficients independent of the parameters. If for some constant $h$ independent of the parameters, and some function $f$ of the parameters, $h(L-\delta)/f$ and $ Q/f^2$ are independently distributed, the former according to the normal law with zero mean and unit variance, the latter according to the $\chi^2$-law with $k-1$ degree of freedom, then the quotient
	$$h(L-\delta)/[Q/(k-1)]^{\frac{1}{2}}$$
	will have the t-distribution  with $k-1$ degree of freedom, no matter what the values of the parameters. And the maximum number of degree of freedom is no greater than $m-1$.
	
\end{theorem}

Thus $T_e$ test have reached the maximum number of degree of freedom.
And in another point of view, it can be understood easily. Even if $m \to +\infty$, or we just know what $\mu_1$ is exactly, the degree of freedom of pivot still depends on $n$, the smaller sample size. It is not the defect of $T_e$ test, but decided by the form of Behrens Fisher Problem.

\subsection*{Division of information}
 Besides, $T_e$ test also show how much information is used and how much information is lost after transformation, which can provide a clearer understanding to this problem:

$$\tilde{X} = \frac{(Q^T_{m\times m})X}{\sqrt{m}} \quad  \tilde{Y} = \frac{(P^T_{n\times n})Y}{\sqrt{n}}$$

After the transformation of orthogonal matrix,
there is no lose of information, but in the final t-test, we use only the first $n$ elements of $\tilde{X}$, and the information of the latter $m-n$ ones have no contributions.

Scheffe also offered a neat form, without giving a clear division of information, though:

\subsection*{Scheffe's transformation}

Suppose now we have two independent normal samples $X_i, i=1,2,\ldots,m$ and $Y_j,j = 1,2,\ldots,n$.  $$X_i \sim N(\mu_1,\sigma^2_1)\quad (i=1,...,m)\quad  Y_j \sim N(\mu_2,\sigma^2_2)\quad (j=1,...,n)$$

Scheffe decided to use the linear combination of $x_i$ and $y_i$ to design an exact confidence interval.
$$d_i = x_i - \sum_{j=1}^{n}c_{ij}y_j \quad (i=1,2,...,m)$$

The variables $d_i$ still follows normal distribution. The necessary and sufficient condition that $d_i$ all share the same mean, equal variance and zero covariance, is that 

$$\sum_{j=1}^{n} c_{ij} = 1,\quad \sum_{k=1}^{n} c_{ik}c_{jk} = c^2 \delta_{ij}$$

Where $\delta_{ii}=1, \delta_{ij}=0$ if $i \not = j$. And Scheffe also give a neat solution

\begin{equation*}
c_{ij}=\left\{
\begin{aligned}
&\delta_{ij}(m/n)^{\frac{1}{2}}-(mn)^{-\frac{1}{2}}+1/n, & j\le m \\
&1/n, & j>m \\
\end{aligned}
\right.
\end{equation*}

Then we can apply a t-test to derive the confidence interval of $d_i$'s mean, which is non-asymptotic and is also the confidence interval of $\mu_1-\mu_2$.\\

%% 有两个版本的证明，因为我目前还证明不了它的凸性（比证明导数为零要难很多），只好先证明一个弱化版本
\section{Is 2-norm best?}

Now that we have reached the maximum number of t-test's degree of freedom, a following problem is that whether t-test with the quadratic form is the best? To be specific about the definition of best, we will use the expectation of confidence interval length to evaluate it.

We will give some definitions at first :

If X follows a normal distribution N(0,1), we will define $|X|^p$'s distribution as $\chi_1^p$; and if $(X_1$, $X_2$, ..., $X_n)$ follows a multivariate normal distribution $N(\vec{0},I_n)$, we will define $\sum_{i=1}^{n} |X_i|^p$'s distribution as $\chi_n^p$. And their probability density functions are as follow.

\begin{equation}
	\label{chi_1}
	f_{\chi^p_1}(x)=\frac{1}{\sqrt{2\pi}}  \frac{2}{p}  exp(-\frac{1}{2}x^{\frac{2}{p}})  x^{\frac{1}{p}-1}
\end{equation}

\begin{equation}
	\label{chi_n}
	f_{\chi^p_n}(y)=\idotsint_{[0,+\infty]^{n-1}}   \frac{1}{(2\pi)^{n/2}  }  \frac{2^n}{p^n}
	exp( -\frac{1}{2} \sum_{i=1}^{n}x_i^{\frac{2}{p} }   )   (\Pi_{i=1}^n x_i)^{ \frac{1}{p}-1} dx_1 dx_2...dx_{n-1}
\end{equation}
where $x_n=y-\sum_{i=1}^{n-1} x_i$.

Now, given normal samples $X_i$, $i=1,2,...,n$, where
$$X_1 \sim N(\mu,\sigma^2) \quad X_i \sim N(0,\sigma^2) \quad (j=2,3,...,n)$$

We can derive a series of $t_p$ test which are defined by
$$t_p: \frac{|X_1-\mu|^p}{\sum_{i=2}^{n} |X_i|^p} \sim  \frac{\chi_1^{p}}{ \chi_{n-1}^{p}  }$$

Then we can derive a confidence interval from it. And we want to find the p that minimize the expectation of confidence interval length.

\begin{theorem}[$t_2$ test reaches the shortest confidence interval, under the assumption of convexity.]
	\label{2min}

	Given a combination of $ n \in Z^+$ and $ \alpha \in (0,1)$,
	for $t_p$ test, we define its corresponding expectation of confidence interval length as l(p):
	
	$$ l(p) =  2q(  (\frac{ \chi_1^p}{\chi_n^p})^{ \frac{1}{p}  } ,\alpha ) E(   ({\chi_n^p})^{ \frac{1}{p}  }  )$$
	
	where
	$$\int_{-\infty}^{q(X,\alpha) }  f_X(x)dx=1-\alpha $$
	$$E(g(X))=\int_{-\infty}^{+\infty} f_X(x)g(x)dx$$

	If we assume the convexity of $l(p)$, then:
	$$l(2)\le l(p), \forall p >0$$
\end{theorem}

\begin{proof}
	
	For a complete proof, we may need to prove that :
	$$\frac{dl}{dp}(2)=0  \quad  \frac{d^2l}{dp^2}(p)\ge0, \forall p>0 $$ 
	%	\begin{figure}[H]
	%		\centering
	%		\subfloat[$n=5$]{\includegraphics[width=0.8\textwidth]{Convex5.jpeg}}
	%		\caption{ 
	%			Above is the simulation of $t_p$ test, which shows the change of %confidence interval length as $p$ changes. 
	%		}
	%		\label{c}
	%	\end{figure}
	Considering the support of simulations, we can just deal with the neighborhood of extreme point and prove:
	$$\frac{dl}{dp}(2)=0  \quad  \frac{d^2l}{dp^2}(2)\ge0$$ 
	What's more, if we have a stronger assumption that the function is convex, all we need to prove is that 
	$$\frac{dl}{dp}(2)=0$$ 
	The proof of this weaker form can be seen in the supplement.
\end{proof}

This theorem points out that given the same data and significance level, $t_2$ test has the shortest confidence interval, in the sense of expectation! In another word, the 2-norm, or the sample variance, contains more information than other p-norm forms when it comes to this problem. And this theorem gives a new strong support to t-test, whose theoretical superiority may be neglected before.
Besides, it is $T_e$ test's clear information division that allow us to derive a $t_p$ test, which also values as one of its contributions.

\section{Simulation studies}

In this section, we compare our method to a few existing methods. We compare their sizes and powers under various settings. To compare the size,  we set $
\mu_1 = \mu_2 = 0$,  and consider the ratios of variances $\sigma_1^2/\sigma_2^2$ ranging over $\{1/25,2/24,3/23,\ldots,24/2,25/1\} $ with the constant $\sigma_1^2+\sigma_2^2=26$.  We set the nominal level $\alpha = 0.05$. These settings are the same as in the paper \citep{ReviewT1}. We consider four pairs of sample sizes with $(m,n) = (50,50)$, $(15,15)$, $(50,5)$ or $(7,5)$. $(m,n) =(50,50)$ means that both data sets have relatively large sample sizes; $(m,n) =(15,15)$ means that both data sets have relatively small sample sizes; $(m,n) =(50,5)$ means that we are studying unbalanced sample sizes; and $(m,n) =(7,5)$ means that we consider unbalanced and extremely small sample sizes.

Type one error and power will serve as the criteria for comparison. We will generate 100,000 samples from two independent normal distributions, with variances and sample sizes not necessarily the same. And when we are simulating the power, we will set $\mu_1-\mu_2=2$.

Besides, a too liberal type one error will make one methods unqualified as a useful test, and thus we will not compare its power to others.
\begin{figure}[H]
	\centering
	\subfloat[Type one Error]{\label{50501}\includegraphics[width=0.5\textwidth]{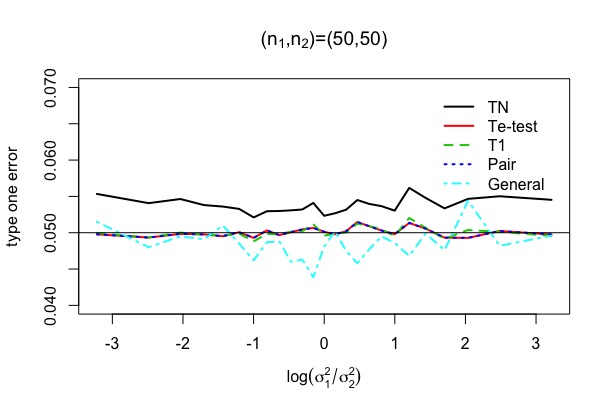}}
	\subfloat[Power]{\label{50502}\includegraphics[width=0.5\textwidth]{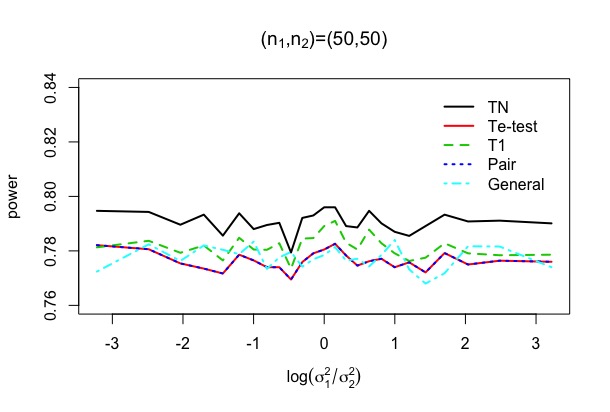}}
	
	\caption{ We can find that when the sample sizes are big, all of five methods perform well. And their powers are quite close. }
	\label{5050}
\end{figure}

\begin{figure}[H]
	\centering
	\subfloat[Type one Error]{\label{15151}\includegraphics[width=0.5\textwidth]{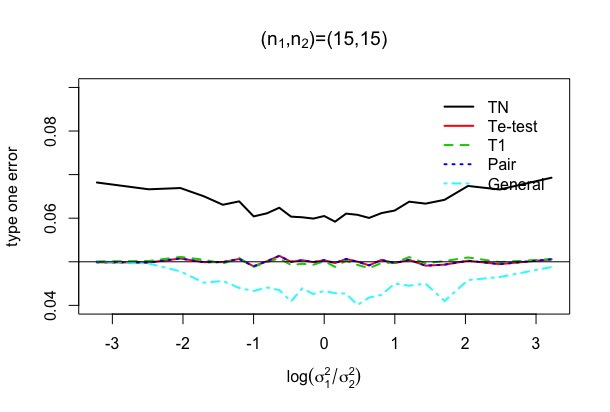}}
	\subfloat[Power]{\label{15152}\includegraphics[width=0.5\textwidth]{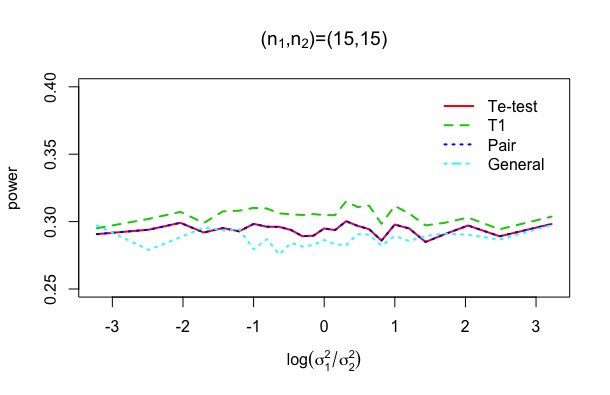}}
	
	\caption{ 
		 We can find that when the sample sizes are smaller, $T_N$ method, which relies on the law of large numbers, performs badly, while the other three's type one errors are still close to $\alpha$. Considering the fact that $T_1$ test is asymptotic, we may get the conclusion that $T_1$ achieve a great fitting to t-distribution with degree of freedom $f$, even when the samples are small(and equal), which is also why it is thought highly of. Besides, its power is even a little bit higher than other non-asymptotic methods', which may result from its degree of freedom $f$ being bigger than $n$. }
	\label{1515}
\end{figure}

\begin{figure}[H]
	\centering
	\subfloat[Type one Error]{\label{15501}\includegraphics[width=0.5\textwidth]{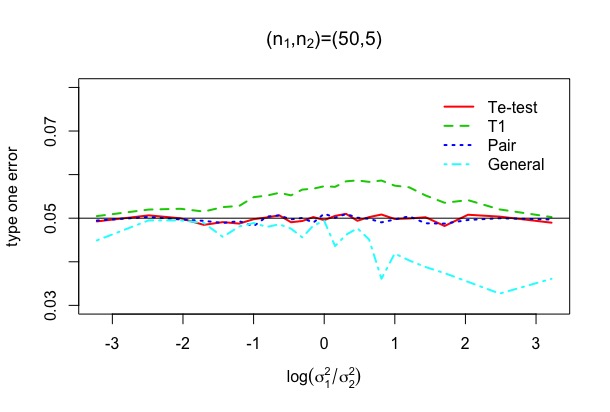}}
	\subfloat[Power]{\label{15502}\includegraphics[width=0.5\textwidth]{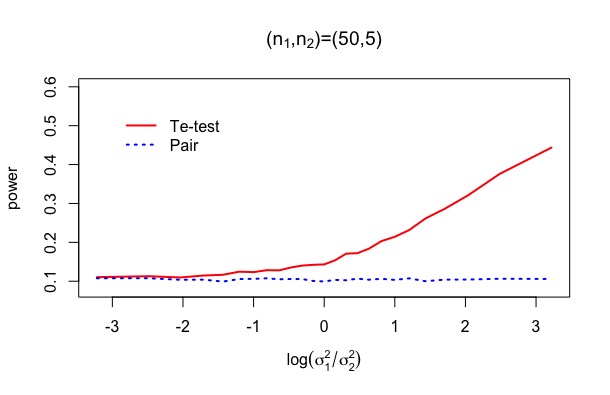}}
	
	\caption{ 
		 We can find that when the imbalance of sample sizes are serious, the $T_1$ test, which is an asymptotic method, fails to fit the model finally. Its type one error is significantly biased from $\alpha$, although its power is still very high. According to Figure \ref{5050} and Figure \ref{1515}, the paired t-test equals to $T_e$ test when the sample sizes are the same, which comes from the proof of $T_e$ test. Now we also find that when the sample sizes are unbalanced, $T_e$ test can use more information of two populations than paired t-test, thus acquiring a higher power, with itself still being a exact test. 
		}
	\label{1550}
\end{figure}

\begin{figure}[H]
	\centering
	\subfloat[Type one Error]{\label{571}\includegraphics[width=0.5\textwidth]{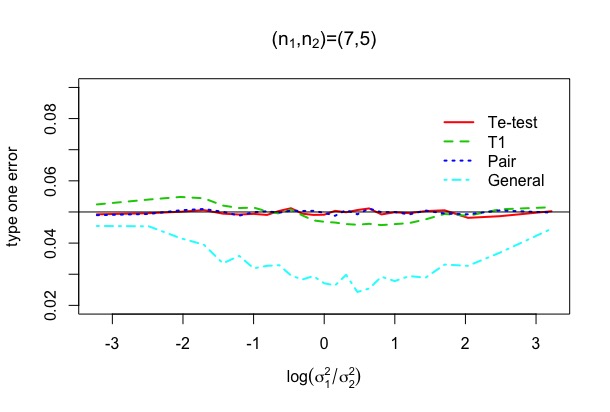}}
	\subfloat[Power]{\label{572}\includegraphics[width=0.5\textwidth]{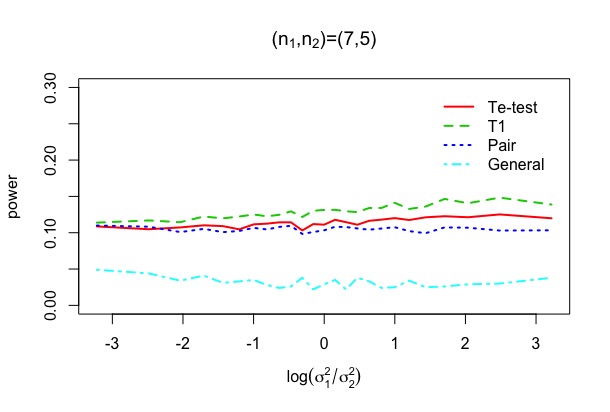}}
	
	\caption{ 
		Even when the sample sizes are extremely small, $T_e$ test still performs well. Although the power is very low, it is not due to test's flaw but results from the limited number of samples.}
	\label{57}
\end{figure}

\section{Development of $T_e$ test}

Here we will look into the condition where samples are from two multi-variate normal distributions, and the variance matrices are not necessarily the same.

$$X_i\sim N(\vec{\mu}_1,\Sigma_1),(i=1,2,...,m)$$
$$Y_i\sim N(\vec{\mu}_2,\Sigma_2),(i=1,2,...,n)$$
where $\Sigma_1$,$\Sigma_2$ are p by p matrices.
$P_{n\times n}$, $Q_{m\times n}$ is the same as those in Te-test.
$$\tilde{X}_{n \times p}=P_{n\times m}X_{m\times p}/\sqrt{m}=
(\tilde{x}_1,\tilde{x}_2,...,\tilde{x}_n)^T$$

$$\tilde{Y}_{n \times p}=Q_{n\times m}Y_{m\times p}/\sqrt{m}=
(\tilde{y}_1,\tilde{y}_2,...,\tilde{y}_n)^T$$

Then we can get two series of variables with independently identically distribution
$$\tilde{x}_1-\vec{\mu}_1,\tilde{x}_2,...,\tilde{x}_n \sim N(0,\Sigma_1/m)$$
$$\tilde{y}_1-\vec{\mu}_2,\tilde{y}_2,...,\tilde{y}_n \sim N(0,\Sigma_2/n)$$
Let $z_1=(\tilde{x}_1-\vec{\mu}_1)-(\tilde{y}_1-\vec{\mu}_2)$,$z_i=\tilde{x}_i-\tilde{y}_i,(i=2,3,...,n)$

According to Wishart distribution, we have
$$S=\sum_{i=2}^{n} z_iz_i^T  \sim W_p(n-1,\frac{\Sigma_1}{m}+\frac{\Sigma_2}{n})$$
$$z_1\sim N(0,\frac{\Sigma_1}{m}+\frac{\Sigma_2}{n})$$

$$z_1   \bot S$$
Then we are able to derive the final pivot in the form of Hotelling's T-Square:

$$(n-1)z_1^TS^{-1}z_1\sim T_{p,n-1}^2 =\frac{p(n-1)}{n-p}F_{p,n-p}$$

We can also get the confidence region of $\vec{\mu}_1-\vec{\mu}_2$ from 
$$[(\tilde{x}_1-\tilde{y}_1)-(\vec{\mu}_1-\vec{\mu}_2)]^TS^{-1}[(\tilde{x}_1-\tilde{y}_1)-(\mu_1-\mu_2)]\sim \frac{p}{n-p}F_{p,n-p}$$

%%%考虑放到 supplement里面

\section{Two stage test}

For methods above, we can obtain no more samples except the given ones. But in some conditions, we may perform the second stage experiment to get more samples. In another word, for any given power or given length of confidence interval, we want to design a two stage test whose level is exactly $\alpha$. 

In the form of experiment design, there are already many kinds of exact tests \citep{twostage1,twostage2,twostage3}. Because their similar performance when variance ratio is close and Chapman's neat form, we will take it into account and compare it to Te-test.

\subsection*{Chapman's: Stage One}
Take initial samples $X_1$, $X_2$, ..., $X_{n_0}$ and $Y_1$, $Y_2$ ,..., $Y_{n_0}$ (both of size $n_0 \ge 2$) from the normal distributions $N(\mu_1,\sigma_1^2)$ and $N(\mu_2,\sigma_2^2)$ respectively, and calculate

$$\bar{X}(n_0)=\frac{1}{n_0}\sum_{i=1}^{n_0}X_i,\quad  \bar{Y}(n_0)=\frac{1}{n_0}\sum_{i=1}^{n_0}Y_i$$

$$S_1^2=\frac{1}{n_0-1}\sum_{i=1}^{n_0}(X_i-\bar{X}(n_0))^2,\quad  S_2^2=\frac{1}{n_0-1}\sum_{i=1}^{n_0}(Y_i-\bar{Y}(n_0))^2$$

$$n_1=\max (n_0+1,[S_1^2/h^2]), \quad n_2=\max (n_0+1,[S_2^2/h^2])$$

where $[x]$ means the smallest integer that is greater than x and h is a given constant(we will refer to more about it in the second stage).

\subsection*{Chapman's: Stage Two}
Take additional $n_1-n_0$ samples $X_{n_0+1}, X_{n_0+2},...,X_{n_1}$from distribution $N(\mu_1,\sigma_1^2)$, and take additional $n_2-n_0$ samples $Y_{n_0+1}, Y_{n_0+2},...,Y_{n_1}$from distribution $N(\mu_2,\sigma_2^2)$. And then calculate

$$\tilde{\bar{X}}=a_1  \sum_{i=1}^{n_0}X_i + a_2  \sum_{i=n_0+1}^{n_1}X_i,\quad 
\tilde{\bar{Y}}=b_1  \sum_{i=1}^{n_0}Y_i + b_2  \sum_{i=n_0+1}^{n_2}Y_i $$

where 

$$  n_0 a_1+(n_1-n_0)a_2=1, \quad n_0 a_1^2+(n_1-n_0)a_2^2=\frac{h^2}{S_2^2}$$
$$  n_0 b_1+(n_1-n_0)b_2=1, \quad n_0 b_1^2+(n_1-n_0)b_2^2=\frac{h^2}{S_2^2}$$

we can prove the existence of the solution $a_i,b_i\ (i=1,2)$ because of the fact $n_i \ge S_i^2/h^2$  .

Then we can prove that $(\tilde{\bar{X}}-\mu_1)/h $ and $(\tilde{\bar{Y}}-\mu_2)/h $ are both Student's-t random variables with $n_0-1$ degrees of freedom respectively. When it comes to the confidence interval of $\mu_1-\mu_2$, we can derive it from 
$$(\tilde{\bar{X}}-\tilde{\bar{Y}}-(\mu_1-\mu_2))/h  \sim t_{n_0-1}-t_{n_0-1}$$
$$\mu_1-\mu_2\in [\tilde{\bar{X}}-\tilde{\bar{Y}}-c_{1-\alpha/2}h, \tilde{\bar{X}}-\tilde{\bar{Y}}+c_{1-\alpha/2}h]$$

where $\tilde{f}$ is the probability density function of the difference between two Student’s-t random variable with $n_0-1$ degrees of freedom and 
$$c_{1-\alpha/2}=\int_{-\infty}^{1-\alpha/2} \tilde{f}(x)dx$$

What's more, when we want to perform a test with condifence interval of given length 2d, we can achieve it easily by setting $h=d/c_{1-\alpha/2}$ at the first stage.

\subsection*{Disadvantage}

When we compare it to $T_e$ test, we can find some shortcomings of it: 
they are in the form of experiment design, so it fails to deal with fixed size questions(i.e. when samples are given); 
the variance information of data in stage two is lost in Chapman's method; 
the decision of $n_0$ depends on specific problem, and improper choice of $n_0$ may make $n_1-n_0$ extremely big, leading to great loss of information. 

\begin{figure}[H]
	\centering
	\subfloat[$n_0=10$]{\includegraphics[width=0.5\textwidth]{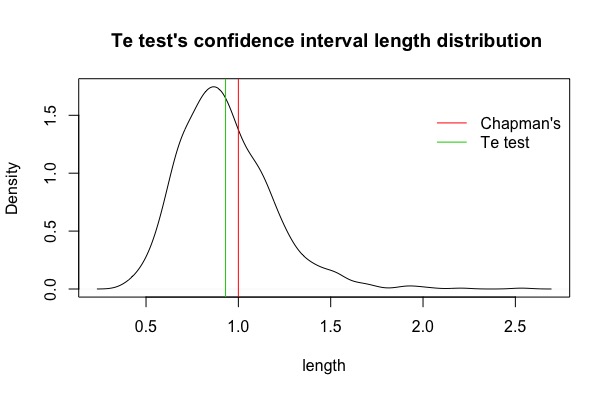}}
	\subfloat[$n_0=20$]{\includegraphics[width=0.5\textwidth]{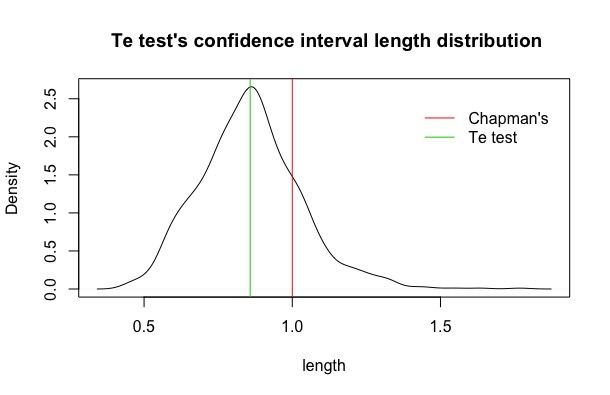}}
	\caption{ 
		Here we set $\mu_1=\mu_2=0$, $\sigma_1^2=1$, $\sigma^2_2=25$. The given confidence interval length of Chapman's method is set to $d=1$ (the red line). We can find that for the same data, the expectation of $T_e$ test's confidence interval length (the green one) is smaller.
	}
	\label{chapman}
\end{figure}

In fact, if we apply Te-test to the samples gained in stage one and two, the expectation of confidence length is usually smaller than that of Chapman's method, which means Te-test can make more use of information (see Figure \ref{chapman}).

\subsection*{Advantage}

For any given confidence interval length, we can design a corresponding exact method. And to the date only two-stage test can achieve it.
Although Te-test can use more information sometimes, we can't perform Chapman's method first and then apply Te-test to the samples gained. Because the sample size depends on the samples collected in stage one, this combined method is no longer exact. Thus two-stage test is still indispensable when it comes to experiment designing.

\section{Conclusion}

Above all, $T_1$ test, which is thought high of, does perform well in Behrens-Fisher Problem when the sample sizes are big, with its power even higher than our $T_e$ test. But its nature of an asymptotic test makes it fail to fit the model when the sample sizes are seriously unbalanced. As for the $T_e$ test, its character of non-asymptotic makes its type one error reach $\alpha$ exactly. Its ability to make more use of all data from two populations contributes to its gaining a higher power than paired t-test's, which is also a non-asymptotic method. And we can see from the simulations that it is robust even when the sample sizes are extremely small.

Of course there are some disadvantages with it. Traditional theory argues that paired t-test makes sense only when two populations are paired. But our $T_e$ test share the same idea with paired t-test, and it is sometimes doubted whether the two populations in Behrens-Fisher Problem are independent. Besides, if we choose different permutations of $x_i$ and $y_j$, we may get different outcomes from $T_e$ test, they are identical distributions though, which may be controversial for a test.

\bibliographystyle{apalike}
\bibliography{paper-ref}

\section{Supplement}

\subsection{Lemma part}

\begin{lemma}[$F^p$'s properties]\label{Fp}

	\begin{equation}
	f_{F^p}(z)= \int_{0}^{+\infty} yf_{\chi_n^p}(y)f_{\chi_1^p}(yz) dy
	\end{equation}
	
	Especially, when $p=2$, $f_{F^p}(z)$ has an explicit expression:
	$$f_{F^2}(z)=\frac{1}{B(\frac{1}{2},\frac{n}{2})} x^{-\frac{1}{2}} (1+x)^{ -\frac{n+1}{2}} $$
	
	\begin{align}
	\frac{d f_{F^p}(z)}{d p} =
	\label{ABC left}&\frac{1}{\sqrt{2\pi}}z^{-\frac{1}{2}}(1+z)^{-\frac{n+1}{2}}
	[ \frac{A  }{p^2(1+z)  } + \frac{ B }{p^2  }\frac{  zlnz}{1+z  }-\frac{  C}{p^2  } +\frac{  D}{p^2  }\frac{  z}{1+z  }-\frac{ E }{p^2  }-\frac{  F}{ p^2 }lnz-\frac{  n+1}{ p }F         ]
	\end{align}

	where
	$$A=\idotsint_{[0,+\infty]^{n}}  (\sum_{i=1}^{n} y_i)^{\frac{1}{2}}
	\frac{1}{(2\pi)^{n/2}  } 
	exp( -\frac{1}{2} \sum_{i=1}^{n}y_i   )   (\Pi_{i=1}^n y_i)^{ \frac{1}{p}-1}  (\sum_{i=1}^{n}y_ilny_i )  dy_1 dy_2...dy_{n-1}dy_n$$
	
	$$B=\idotsint_{[0,+\infty]^{n}}  (\sum_{i=1}^{n} y_i)^{\frac{1}{2}}
	\frac{1}{(2\pi)^{n/2}  } 
	exp( -\frac{1}{2} \sum_{i=1}^{n}y_i   )   (\Pi_{i=1}^n y_i)^{ \frac{1}{p}-1}  (\sum_{i=1}^{n}y_i )  dy_1 dy_2...dy_{n-1}dy_n$$
	
	$$C=\idotsint_{[0,+\infty]^{n}}  (\sum_{i=1}^{n} y_i)^{\frac{1}{2}}
	\frac{1}{(2\pi)^{n/2}  } 
	exp( -\frac{1}{2} \sum_{i=1}^{n}y_i   )   (\Pi_{i=1}^n y_i)^{ \frac{1}{p}-1}  (\sum_{i=1}^{n}lny_i )  dy_1 dy_2...dy_{n-1}dy_n$$
	
	$$D=\idotsint_{[0,+\infty]^{n}}  (\sum_{i=1}^{n} y_i)^{\frac{1}{2}}
	\frac{1}{(2\pi)^{n/2}  } 
	exp( -\frac{1}{2} \sum_{i=1}^{n}y_i   )   (\Pi_{i=1}^n y_i)^{ \frac{1}{p}-1}  (\sum_{i=1}^{n}y_i )ln(\sum_{i=1}^{n} y_i )  dy_1 dy_2...dy_{n-1}dy_n$$
	
	$$G=\idotsint_{[0,+\infty]^{n}}  (\sum_{i=1}^{n} y_i)^{\frac{1}{2}}
	\frac{1}{(2\pi)^{n/2}  } 
	exp( -\frac{1}{2} \sum_{i=1}^{n}y_i   )   (\Pi_{i=1}^n y_i)^{ \frac{1}{p}-1}  ln(\sum_{i=1}^{n}y_i)   dy_1 dy_2...dy_{n-1}dy_n$$
	
	$$F=\idotsint_{[0,+\infty]^{n}}  (\sum_{i=1}^{n} y_i)^{\frac{1}{2}}
	\frac{1}{(2\pi)^{n/2}  } 
	exp( -\frac{1}{2} \sum_{i=1}^{n}y_i   )   (\Pi_{i=1}^n y_i)^{ \frac{1}{p}-1}    dy_1 dy_2...dy_{n-1}dy_n$$
\end{lemma}

\begin{proof}

\begin{align}
	 \frac{df_{F^p}(z)}{dp} &= \int_{0}^{+\infty} y [ \frac{d f_{\chi_n^p}(y)}{dp}f_{\chi_1^p}(yz)        +f_{\chi_n^p}(y) \frac{df_{\chi_1^p}(yz)}{dp}    ] dy\\
	=&         \label{transform}
	\idotsint_{[0,+\infty]^{n}}  (\sum_{i=1}^{n} x_i)
	\frac{1}{(2\pi)^{n/2}  }  \frac{2^n}{p^n}
	exp( -\frac{1}{2} \sum_{i=1}^{n}x_i^{\frac{2}{p} }   )   (\Pi_{i=1}^n x_i)^{ \frac{1}{p}-1}\nonumber \\
	& \frac{1}{\sqrt{2\pi}}  \frac{2}{p}  exp(-\frac{1}{2}(z\sum_{i=1}^{n} x_i)^{\frac{2}{p}})  (z\sum_{i=1}^{n} x_i)^{\frac{1}{p}-1} \nonumber \\
	&[\frac{1}{p^2}\sum_{i=1}^{n} (x_i^{ \frac{2}{p}  }-1)lnx_i -\frac{n}{p}  + 
	\frac{1}{p^2}(z\sum_{i=1}^{n} x_i -1)ln(z\sum_{i=1}^{n} x_i) -\frac{1}{p} ]dx_1 dx_2...dx_{n-1}dx_n\\
	=& \label{1+z}
	\idotsint_{[0,+\infty]^{n}}  (\sum_{i=1}^{n} y_i)^{\frac{1}{2}}
	\frac{1}{(2\pi)^{(n+1)/2}  } 
	exp( -\frac{1}{2} \sum_{i=1}^{n}y_i   )   (\Pi_{i=1}^n y_i)^{ \frac{1}{p}-1}  z^{-\frac{1}{2}}(1+z)^{-\frac{n+1}{2}}  \nonumber  \\
	&\frac{1}{p^2}  [\sum_{i=1}^{n} (\frac{y_i}{1+z}-1)ln\frac{y_i}{1+z} -pn  + 
	(\frac{z}{1+z}\sum_{i=1}^{n} y_i -1)ln(\frac{z}{1+z}\sum_{i=1}^{n} y_i) -p ]dy_1 dy_2...dy_{n-1}dy_n\\
	=&  \label{ABC}
	\frac{1}{\sqrt{2\pi}}z^{-\frac{1}{2}}(1+z)^{-\frac{n+1}{2}} [\frac{1}{p^2(1+z)}A-\frac{1}{p^2}(\frac{ln(1+z)}{1+z}-\frac{z}{1+z}ln\frac{z}{1+z})B-\frac{1}{p^2}C+\frac{1}{p^2}\frac{z}{1+z}D \nonumber \\
	&-\frac{1}{p^2}G+
	(\frac{nln(1+z)}{p^2}-\frac{ln(z)}{p^2}+\frac{ln(1+z)}{p^2}-\frac{n+1}{p})F]\\
	=&\frac{1}{\sqrt{2\pi}}z^{-\frac{1}{2}}(1+z)^{-\frac{n+1}{2}}
	[ \frac{A  }{p^2(1+z)  } + \frac{ B }{p^2  }\frac{  zlnz}{1+z  }-\frac{  C}{p^2  } +\frac{  D}{p^2  }\frac{  z}{1+z  }-\frac{ E }{p^2  }-\frac{  F}{ p^2 }lnz-\frac{  n+1}{ p }F         ]\\
\end{align}

In equation \ref{transform}, we  transform the variates $(x_1,x_2,...,x_{n-1},y)$ to $(x_1,x_2,...,x_n)$; in equation \ref{1+z}, we use $x_i(1+z)=y_i$ to simplify the result.

\end{proof}

\begin{lemma}[$\chi_n^p$'s properties]
	when $p=2$, $$E[   ({\chi_n^p})^{ \frac{1}{p}  }  ]= \int_{0}^{+\infty}  x^{\frac{1}{2}}     \frac{x^{\frac{n}{2}-1 } e^{-\frac{x}{2}}   }    {2^{\frac{n}{2}}\Gamma(\frac{n}{2})} dx = \frac{\sqrt{2\pi}}{B(\frac{1}{2},\frac{n}{2})}=F$$
	
	$$\frac{dE[   ({\chi_n^p})^{ \frac{1}{p}  }  ]}{dp}=\int_{0}^{+\infty} \frac{df_{\chi_n^p}(y)}{dp} y^{\frac{1}{p}} + f_{\chi_n^p}(y)(lny)(-\frac{1}{p^2})y^{-\frac{1}{p}} dy = \frac{1}{p^2} A - \frac{1}{p^2} C -\frac{n}{p}F - \frac{1}{p^2}G$$
\end{lemma}

\begin{proof}
	It is similar to the proof of lemma \ref{Fp}, and we skip it for simplicity.
\end{proof}
\begin{lemma}\label{BF}
	$$B=(n+1)F$$
	
\end{lemma}
\begin{proof}
	
	\begin{align*}
	F&=\idotsint_{[0,+\infty]^{n}}  1\times (\sum_{i=1}^{n} y_i)^{\frac{1}{2}}
	\frac{1}{(2\pi)^{n/2}  } 
	exp( -\frac{1}{2} \sum_{i=1}^{n}y_i   )   (\Pi_{i=1}^n y_i)^{ \frac{1}{p}-1}  dy_1 dy_2...dy_{n-1}dy_n\\
	&=\int_{0}^{+\infty}  x^{\frac{1}{2}}     \frac{x^{\frac{n}{2}-1 } e^{-\frac{x}{2}}   }    {2^{\frac{n}{2}}\Gamma(\frac{n}{2})} dx \\
	&=\frac{2^{\frac{1}{2}}\Gamma(\frac{n+1}{2})}{\Gamma(\frac{n}{2})}=\frac{\sqrt{2\pi}}{   B(\frac{1}{2},\frac{n}{2}  )     }
	\end{align*}
	
	\begin{align*}
	B&=\idotsint_{[0,+\infty]^{n}}  (\sum_{i=1}^{n} y_i)^{\frac{1}{2}}
	\frac{1}{(2\pi)^{n/2}  } 
	exp( -\frac{1}{2} \sum_{i=1}^{n}y_i   )   (\Pi_{i=1}^n y_i)^{ \frac{1}{p}-1}  (\sum_{i=1}^{n}y_i )  dy_1 dy_2...dy_{n-1}dy_n\\
	&=\int_{0}^{+\infty} x \times x^{\frac{1}{2}}     \frac{x^{\frac{n}{2}-1 } e^{-\frac{x}{2}}   }    {2^{\frac{n}{2}}\Gamma(\frac{n}{2})} dx \\
	&=\frac{2^{\frac{3}{2}}\Gamma(\frac{n+3}{2})}{\Gamma(\frac{n}{2})}=(n+1)F
	\end{align*}
\end{proof}

\begin{lemma}\label{DG}
	
	$$D=(n+1)G+2F$$
	
\end{lemma}
\begin{proof}
	
	\begin{align*}
	D&=\idotsint_{[0,+\infty]^{n}}  (\sum_{i=1}^{n} y_i)^{\frac{1}{2}}
	\frac{1}{(2\pi)^{n/2}  } 
	exp( -\frac{1}{2} \sum_{i=1}^{n}y_i   )   (\Pi_{i=1}^n y_i)^{ \frac{1}{p}-1}  (\sum_{i=1}^{n}y_i )ln(\sum_{i=1}^{n} y_i )  dy_1 dy_2...dy_{n-1}dy_n\\
	&=\int_{0}^{+\infty} xlnx \times x^{\frac{1}{2}}     \frac{x^{\frac{n}{2}-1 } e^{-\frac{x}{2}}   }    {2^{\frac{n}{2}}\Gamma(\frac{n}{2})} dx \\
	&=\frac{1 }    {2^{\frac{n}{2}}\Gamma(\frac{n}{2})} \int_{0}^{+\infty}      x^{\frac{n+1}{2} } e^{-\frac{x}{2}} lnx dx \\
	&=\frac{1 }    {2^{\frac{n}{2}}\Gamma(\frac{n}{2})} \int_{0}^{+\infty}      (-2)  x^{\frac{n+1}{2} }  lnxde^{-\frac{x}{2}} \\
	&=\frac{1 }    {2^{\frac{n}{2}}\Gamma(\frac{n}{2})} \int_{0}^{+\infty}      (n+1)x^{\frac{n-1}{2} } e^{-\frac{x}{2}} lnx + 2x^{\frac{n-1}{2} } e^{-\frac{x}{2}}  dx \\
	&=(n+1)G+2F
	\end{align*}
\end{proof}

\begin{lemma}\label{L2}
	$$\int_{0}^{1}  x^{u-1}(1-x)^{v-1} lnx dx = B(u,v)[\psi(u)-\psi(u+v)]$$
	where $$\psi(z) = -c + \sum_{n=0}^{+\infty} (\frac{1}{n+1} - \frac{1}{z+n})$$
\end{lemma}

\begin{lemma}\label{L1}
	$$\int_{0}^{+\infty}  x^{v-1} exp(-ux)lnx dx = \frac{1}{u^v}\Gamma(v)[\psi(v)-lnu)] $$
\end{lemma}

\begin{lemma}\label{AC}
	$$nA=(n+1)C+(2n^2+2n-2)F$$
\end{lemma}
\begin{proof}

	\begin{align*}
	\frac{A}{n}&=\idotsint_{[0,+\infty]^{n}}  (\sum_{i=1}^{n} y_i)^{\frac{1}{2}}
	\frac{1}{(2\pi)^{n/2}  } 
	exp( -\frac{1}{2} \sum_{i=1}^{n}y_i   )   (\Pi_{i=1}^n y_i)^{ \frac{1}{p}-1}  (y_1lny_1)  dy_1 dy_2...dy_{n-1}dy_n\\
	&=\int_{0}^{+\infty} \int_{0}^{+\infty}  \frac{1}{\sqrt{2\pi}} (x+y)^{\frac{1}{2}}  \frac{y^{\frac{n-1}{2}-1 } e^{-\frac{y}{2}}   }    {2^{\frac{n-1}{2}}\Gamma(\frac{n-1}{2})} e^{-\frac{1}{2}x} x^{-\frac{1}{2}}xlnx \ dydx\\
	&= \int_{0}^{+\infty} \int_{0}^{t}   \frac{1}{\sqrt{2\pi}}  \frac{1}{2^{\frac{n-1}{2}}\Gamma(\frac{n-1}{2})} e^{-\frac{1}{2}t}t^{\frac{1}{2}} y^{\frac{n-3}{2}}(t-y)^{\frac{1}{2}}ln(t-y)\ dydt\\
	&=\int_{0}^{+\infty} \int_{0}^{1}   \frac{1}{\sqrt{2\pi}}  \frac{1}{2^{\frac{n-1}{2}}\Gamma(\frac{n-1}{2})} e^{-\frac{1}{2}t}t^{\frac{1}{2}}t^{\frac{n-3}{2}} (1-x)^{\frac{n-3}{2}}t^{\frac{1}{2}}x^{\frac{1}{2}}(lnx+lnt)\ tdxdt\\
	&=\frac{1}{\sqrt{2\pi}}  \frac{1}{2^{\frac{n-1}{2}}\Gamma(\frac{n-1}{2})}
	\int_{0}^{+\infty}     e^{-\frac{1}{2}t}t^{\frac{n+1}{2}}
	[\int_{0}^{1}   (1-x)^{\frac{n-3}{2}}x^{\frac{1}{2}}lnx \ dx + lnt\int_{0}^{1}   (1-x)^{\frac{n-3}{2}}x^{\frac{1}{2}} \ dx ]
	dt\\
	&=\frac{1}{\sqrt{2\pi}}  \frac{1}{2^{\frac{n-1}{2}}\Gamma(\frac{n-1}{2})}
	\int_{0}^{+\infty}     e^{-\frac{1}{2}t}t^{\frac{n+1}{2}} B(\frac{3}{2},\frac{n-1}{2})[\psi(\frac{3}{2}) - \psi(\frac{n+2}{2})+ lnt] \ dt\\
&=\frac{1}{\sqrt{2\pi}}  \frac{1}{2^{\frac{n-1}{2}}\Gamma(\frac{n-1}{2})}
B(\frac{3}{2},\frac{n-1}{2})[\psi(\frac{3}{2}) - \psi(\frac{n+2}{2})   + \psi(\frac{n+3}{2}) + ln2 ]2^{\frac{n+3}{2}}\Gamma(\frac{n+3}{2})
	\end{align*}
	
	The first four equations are simple variable substitutions;
	the sixth equation comes from Lemma \ref{L1}; the seventh equation comes from Lemma \ref{L2}.
	
	Similarly,
		\begin{align*}
		\frac{C}{n}=&
		\frac{1}{\sqrt{2\pi}}  \frac{1}{2^{\frac{n-1}{2}}\Gamma(\frac{n-1}{2})}B(\frac{1}{2},\frac{n-1}{2})[\psi(\frac{1}{2}) - \psi(\frac{n}{2}) + \psi(\frac{n+1}{2}) + ln2   ]2^{\frac{n+1}{2}}\Gamma(\frac{n+1}{2})\\ 
		\end{align*}
		
 Thus,
 \begin{align*}
 &nA-(n+1)C\\=&
 \frac{  n^2 }{\sqrt{2\pi}} \frac{2^{\frac{n+3}{2}}\Gamma(\frac{n+3}{2})}{2^{\frac{n-1}{2}}\Gamma(\frac{n-1}{2})}B(\frac{3}{2},\frac{n-1}{2})[\psi(\frac{3}{2}) - \psi(\frac{n+2}{2})   + \psi(\frac{n+3}{2}) + ln2 ]\\ 
 -&\frac{  n(n+1) }{\sqrt{2\pi}} \frac{2^{\frac{n+1}{2}}\Gamma(\frac{n+1}{2})}{2^{\frac{n-1}{2}}\Gamma(\frac{n-1}{2})}B(\frac{1}{2},\frac{n-1}{2})[\psi(\frac{1}{2}) - \psi(\frac{n}{2})   + \psi(\frac{n+1}{2}) + ln2 ]\\
 =&	\frac{  n(n+1) }{\sqrt{2\pi}} \frac{2^{\frac{n+1}{2}}\Gamma(\frac{n+1}{2})}{2^{\frac{n-1}{2}}\Gamma(\frac{n-1}{2})} \frac{\Gamma(\frac{1}{2}) \Gamma(\frac{n-1}{2})  }{\Gamma(\frac{n}{2}) }
 [\psi(\frac{3}{2}) - \psi(\frac{n+2}{2})   + \psi(\frac{n+3}{2})  -\psi(\frac{1}{2}) + \psi(\frac{n}{2})   -\psi(\frac{n+1}{2})  ]
 \\
 =& n(n+1)  \frac{2^{\frac{1}{2}}\Gamma(\frac{n+1}{2})}{\Gamma(\frac{n}{2})}[2-\frac{2}{n}+\frac{2}{n+1}] \\
 =&(2n^2+2n-2)\frac{2^{\frac{1}{2}}\Gamma(\frac{n+1}{2})}{\Gamma(\frac{n}{2})}\\
 =&(2n^2+2n-2)F
 \end{align*}

\end{proof}

\subsection{Theorem proof}

\begin{theorem}[$t_2$ test reaches the shortest confidence interval, under the assumption of convexity.]
	\label{2min}
	
	Given a combination of $ n \in Z^+$ and $ \alpha \in (0,1)$,
	for $t_p$ test, we define its corresponding expectation of confidence interval length as l(p):
	
	$$ l(p) =  2q(  (\frac{ \chi_1^p}{\chi_n^p})^{ \frac{1}{p}  } ,\alpha ) E(   ({\chi_n^p})^{ \frac{1}{p}  }  )$$
	
	where
	$$\int_{-\infty}^{q(X,\alpha) }  f_X(x)dx=1-\alpha $$
	$$E(g(X))=\int_{-\infty}^{+\infty} f_X(x)g(x)dx$$

	Then,
	$$\frac{d l(p)}{dp}\Big \rvert_{p=2}=0.$$
	
	And if we have additional validation of convexity that $l''(p)>0$, we have:
	
	$$l(2)\le l(p), \forall p >0.$$
\end{theorem}

\begin{proof}
	
	In the following proof, and we will show that 
	$$\frac{d l(p)}{dp}\Big \rvert_{p=2}=0.$$

	Firstly, according to the definition of quantile:
	\begin{equation*}
	q(  (\frac{ \chi_1^p}{\chi_n^p})^{ \frac{1}{p}  } ,\alpha ) E(   ({\chi_n^p})^{ \frac{1}{p}  }  )= q(  \frac{ \chi_1^p}{\chi_n^p} ,\alpha )^{ \frac{1}{p}  }  E(   ({\chi_n^p})^{ \frac{1}{p}  }  )
	\end{equation*}
	
	Besides, when $f= f_{F_p}$, which is a probability density function, we have:
	$$\int_{0}^{q}f(x,p) dx = 1-\alpha \Rightarrow \int_{0}^{q} \frac{df}{dp} dx +\frac{dq}{dp} f(q,p)=0 \Rightarrow \frac{dq}{dp}=\frac{-\int_{0}^{q} \frac{df}{dp}dx }{f(q,p)}$$
	
	Thus 
	\begin{align}
	\frac{d l}{dp}(2)=0&\Leftrightarrow \frac{ d (q^{\frac{1}{p}} E)}{dp}=0 \Leftrightarrow \frac{d(\frac{lnq}{p}+lnE)}{dp}=0\nonumber \\
	&\Leftrightarrow   \frac{    1   }{qp}  \frac{dq}{dp} - \frac{lnq}{p^2} + \frac{ 1  }{E}\frac{dE}{dp}=0\nonumber \\
	&\Leftrightarrow \frac{-\int_{0}^{q} \frac{df}{dp}dx }{f(q,p)}=\frac{q}{p}lnq-\frac{pq}{E}\frac{dE}{dp}\nonumber \\
	& \label{object1}\Leftrightarrow  \frac{df}{dp}(q,p)=
	\frac{d}{dq}\left( (\frac{pq}{E}\frac{dE}{dp}-\frac{q}{p}lnq)f(q,p)  \right)  
	\end{align}

	When $p=2$, let $a=\frac{p}{E}\frac{dE}{dp}$, $b=-\frac{1}{p}$, 
	
	\begin{align}
	&	\frac{d}{dq}\left( (\frac{pq}{E}\frac{dE}{dp}-\frac{q}{p}lnq)f(q,p)  \right)\nonumber \\
	=&\frac{d}{dq}\left(      \frac{1}{B(\frac{1}{2},\frac{n}{2})} q^{-\frac{1}{2}} (1+q)^{ -\frac{n+1}{2}}   (aq+bqlnq)  \right)\nonumber \\
	=& \label{ABC right}  \frac{1}{B(\frac{1}{2},\frac{n}{2})} q^{-\frac{1}{2}} (1+q)^{ -\frac{n+1}{2}}   [\frac{lnq}{1+q}  (\frac{n+1}{2}b)  + lnq(-\frac{n}{2}b) +\frac{1}{1+q} (\frac{n+1}{2}a)  +(b-\frac{n}{2}a)  ] 
	\end{align}
	
	To prove equation \ref{object1}, we just need to prove that formula \ref{ABC left} is equal to formula \ref{ABC right}, which is equivalent to satisfying the following equation set according to Lemma \ref{Fp}:
	
	\begin{equation}\label{check}
	\left\{
	\begin{aligned}
	\frac{B}{p^2}-\frac{F}{p^2}&=-\frac{n}{2}bF\\
	\frac{D-C-G}{p^2}-\frac{(n+1)F}{p}&=(b-\frac{n}{2}a)F\\
	-\frac{B}{p^2}&=\frac{n+1}{2}bF\\
	\frac{A-D}{p^2}&=\frac{n+1}{2}aF
	\end{aligned} \right.
	\end{equation}
	
	Since $$a=\frac{p}{E}\frac{dE}{dp}=\frac{A-C-G-npF}{pF}$$
	
	Combine it with lemma \ref{BF}, lemma \ref{DG} and lemma \ref{AC}, it is easy to check the validity of equation set (\ref{check}), which finished the proof of this theorem.

\end{proof}

\end{document}